\documentclass[12pt]{article}

\usepackage{amscd,amsthm,amsmath,amssymb,amsfonts}

\newtheorem{Theorem}{Theorem}[section]

\newtheorem{Lemma}[Theorem]{Lemma}

\def\buildreb#1\over#2{\mathrel{\mathop{\kern0pt #1}\limits_{#2}}}

\newtoks\baspagetitre
\baspagetitre={\hfill{}}

\newif\ifpagetitre \pagetitretrue 
\def\footnoterule{\kern -3pt\hrule width 2truein \kern 2.4pt}
\renewcommand{\thepage}{\ifpagetitre\the\baspagetitre\global\pagetitrefalse
 \else \arabic{page}\fi}
 \def\eps{\varepsilon} \def\alp{\alpha} \def\Z{{\mathbb Z}} \def\vph{\varphi} \def\R{{\mathbb R}}
 \def\Q{{\mathbb Q}} \def\gam{\gamma} \def\sig{\sigma}   \def\C{{\mathbb C}}
 \def\F{{\mathbb F}}

\begin{document}
  \title{\sc Simultaneous multiplicative rational approximation
 to a real   and a   $p$-adic numbers}  
   
  \author{\bf Yann Bugeaud \& Bernard de Mathan}
 
  \date{}
\maketitle
{\vskip 15mm}

\leftskip = 0pt
\rightskip = 2pt
{\bf R\'esum\'e.} Nous donnons de nouveaux exemples de paires form\'ees d'un nombre r\'eel et d'un nombre $p$-adique  
qui v\'erifient une conjecture d'Einsiedler et Kleinbock sur l'approximation multiplicative simultan\'ee.   
\par {\bf Abstract.}
We give new examples of pairs composed of a real and a $p$-adic numbers that satisfy  a conjecture  on  
simultaneous multiplicative approximation by rational numbers formulated by  Einsiedler and Kleinbock in 2007.  

\par
\ \par
{\bf Keywords} Littlewood conjecture. Simultaneous diophantine approximation. Non-archimedean valuations.\par

\section{Introduction} Let $p$ be a prime number and $\Q_p$ denote the $p$-adic field. 
The $p$-adic absolute value is normalized by $|p|_p=1/p$. It is well known that the diagonal of $\Q\times\Q$ is everywhere dense in $\R\times\Q_p$. The conjecture of Einsiedler and Kleinbock \cite{Eis}  is the following: every nonzero pair $(\alp,\beta)$ in $\R\times\Q_p$ 
satisfies    
$$
\inf_{\substack{q, r \in \Z \\ qr\ne0}} |q||q\alp-r||q\beta-r|_p=0. 
\eqno\hbox{(EK)}  
$$   
Einsiedler and Kleinbock \cite{Eis} proved that the set of exceptions to (EK) has Hausdorff dimension 0; 
however, their conjecture is still open.  
A stronger question  
can be considered, namely, is it true that every nonzero pair $(\alp,\beta)$ in $\R\times\Q_p$ satisfies  
$$\inf_{\substack{q, r \in \Z \\ qr\ne0}}  \max\{|q|, |r|\} |q\alp-r||q\beta-r|_p=0\ ?\eqno\hbox{(EK+)}$$
Badziahin and Bugeaud \cite{BuBa} noticed  that  (EK+) is {\sl not} satisfied when there exists an irreducible polynomial $f(X)$ in $\Z[X]$ of degree 2 such that $f(\alp)=0$ and $f(\beta)=0$. For convenience, we reproduce their argument.   
Let $\alp'$ be the conjugate of $\alp$ in $\R$ and $\beta'$
that of $\beta$ in $\Q_p$. In $\R[X]$, write 
$$f(X)=a(X-\alp)(X-\alp'),$$ where $a$ is a nonzero integer. In $\Q_p[X]$, we have:
$$f(X)=a(X-\beta)(X-\beta').$$ Then, given a pair of integers $(q,r)$ with $q>0$, we have
$$
0<q^2|f(\frac rq)|\le|a|(q|\alp'|+|r|)|q\alp-r|\le|a|(|\alp'|+1)\max\{q,|r|\}|q\alp-r|, \eqno\hbox{(1.1)}
$$ while
$$|q^2f(\frac rq)|_p\le|a|_p\max\{1,|\beta'|_p\}|q\beta-r|_p.\eqno\hbox{(1.2)}$$ 
As $q^2f(r/q)$ is a nonzero integer, we have 
$$1\le|q^2f(\frac rq)| |q^2f(\frac rq)|_p\ ,$$ thus we deduce from (1.1) and (1.2) that
$$\max\{q,|r|\}|q\alp-r||q\beta-r|_p\ge{1\over|a||a|_p(|\alp'|+1)\max\{1,|\beta'|_p\}}\cdot$$ 
We do not know whether there exist other counterexamples to (EK+)  
than those given above;    
however, it is proved in \cite{Bu} that the set of counterexamples to (EK+) has Hausdorff dimension $0$.

Some particular cases have been studied. For instance, if $\alp$ is a nonzero real number,  
then the pair $(1/\alp,0)$ satisfies (EK) if and only if $\alp$ satisfies the {\sl mixed} Littlewood conjecture
formulated by de Mathan and Teuli\'e \cite{dMT}, that is, if and only if,   
$$\inf_{\substack{q, r \in \Z \\ r\ne0}}    |r||r|_p|r\alp-q|=0. \eqno\hbox{(M)}$$  Indeed, if $0<|r|\le|q|/(2|\alp|)$, 
then we have $|q||q-r\alp||r|_p\ge|\alp||r||r|_p\ge|\alp|$, thus if
$$\inf_{qr\ne0}|q||q-r\alp||r|_p=0,$$ we can suppose in the definition of the infimum that $|r|>|q|/(2|\alp|)>0$, and then we get (M); conversely, if we have (M), we also have  
$$\inf_{r\ne0,|q|\le2|r\alp|}|r||r|_p|r\alp-q|=0,$$ which ensures (EK+) for
 $(1/\alp,0)$. 
 Note that in the case of a pair $(\alp,0)$,    
 (EK) and (EK+) are equivalent. It is proved in \cite{dMT} that if $\alp$ is a real quadratic number, then the mixed conjecture is satisfied, thus the pair $(\alp,0)$ satisfies (EK+).
Note also that it has been proved in \cite{Adal} that an analogue of the mixed Littlewood conjecture fails in the field of formal power series $\F_3((T^{-1}))$.\par 
 If we consider the pair $(0,\beta)$, where $\beta\in\Q_p^\star$, then we get the {\sl multiplicative} conjecture,
 studied by de Mathan \cite{BM1,Mon}:
 $$\inf_{qr\ne0}|qr||q\beta-r|_p=0.\eqno\hbox{(MU)}$$ Obviously $\beta$ satisfies (MU) if and only if $1/\beta$ satisfies (MU).
 Let us consider the conditions $$\inf_{0<|q|\le|r|}|qr||q\beta-r|_p=0\eqno\hbox{(MU1)}$$ and
  $$\inf_{0<|r|\le|q|}|qr||q\beta-r|_p=0.\eqno\hbox{(MU2)}$$
 Then (MU) is satisfied if and only if one at least of the conditions (MU1) and (MU2) is satisfied, and a nonzero $p$-adic number $\beta$ satisfies (MU2) if and only if $1/\beta$ satisfies (MU1).\par
 Condition (MU2) implies (EK+) for $(0,\beta)$. 
We do not know whether condition
(MU) implies (MU2). 
\par
We notice that, when $\beta$ satisfies (MU1), then for each real number $\alp$, condition (EK) is satisfied, since $|q\alp-r|\le(|\alp|+1)|r|$.\par 
It has been proved in \cite{BM1} that if $\beta$ is a quadratic number in $\Q_p$, then there exist a positive real number $C$ and infinitely many pairs of integers $(q,r)$ with 
 $$0<|r|\le C{|q|\over\log\max\{|q|,2\}} ,\qquad|q\beta-r|_p\le C|q|^{-2},$$ which implies that
 $$\max\{|q|,|r|\}|r||q\beta-r|_p\le {C^2\over\log\max\{|q|,2\}}\cdot$$ 
 Thus $\beta$ and $1/\beta$ satisfy (MU1), that is to say that $\beta$ satisfies both (MU1) and (MU2). 
 Consequently, the pair $(\alp,\beta)$ satisfies (EK) for each $\alp\in\R$. 
 However, as we have seen above, if $\alpha$ is a real root of the minimal defining polynomial of $\beta$, then   
 the pair $(\alp,\beta)$ does not satisfy (EK+).    
 
 We use the Vinogradov notations $\ll$, $\gg$ and $\asymp$ (we write $A\ll B$ to mean that $0\le A\le KB$ for some positive constant $K$; $A\asymp B$ means that $A\ll B$ and $B\ll A$). The dependence of the implied constants will be generally clear, sometimes 
 we indicate it.    
 Our main theorems are the following extensions of the results of \cite{BM1,dMT}.    
 
\begin{Theorem} Let $\beta$ be a quadratic number in $\Q_p$ and $f(X)$ an irreducible polynomial in $\Z[X]$ such that $f(\beta)=0$. 
Then, for each real number $\alp$ with $f(\alp)\ne0$, the pair $(\alp,\beta)$ satisfies (EK+). 
More precisely, we have 
 $$\buildreb{\lim\inf}\over{\max\{|Q|,|R| \} \to + \infty}    
 \max\{|Q|,|R|\}  \log\max\{|Q|,|R|\}  |Q\alp-R||Q\beta-R|_p <+\infty.$$ 
 \end{Theorem} 
 
 Actually, the proof of Theorem 1.1 shows that there exist  
 infinitely many pairs of integers $(Q,R)$ with   
 $$|Q\alp-R|\log\max\{|Q|,|R|\}\ll_{\alp,\beta}\max\{|Q|,|R|\},\qquad|Q\beta-R|_p\ll_{\alp,\beta} \max\{|Q|,|R|\}^{-2}.\eqno\hbox{(1.3)}$$   
 For $\alp=0$, we obtain for $|Q|>1$, 
 $$|R|\ll{|Q|\over \log|Q|},\qquad|Q\beta-R|_p\ll|Q|^{-2},$$  
 which was proved in \cite{BM1}.\par
 A ``dual'' result is also true.
 
\begin{Theorem}
Let $\alp$ be a quadratic real number and $f(X)$ an irreducible polynomial in $\Z[X]$ such that $f(\alp)=0$. Let $\beta$ be a $p$-adic number.  
If $f(\beta)\ne0$, then the pair $(\alp,\beta)$ satisfies (EK+). 
More precisely, we have
$$\buildreb{\lim\inf}\over{\max\{|Q|,|R|\} \to +\infty}
\max\{|Q|,|R|\}\log\max\{|Q|,|R|\}|Q\alp-R||Q\beta-R|_p<+\infty.$$ 
\end{Theorem}

Actually, the proof of Theorem 1.2 
shows that there exist   
 infinitely many pairs of integers $(Q,R)$ such that
$$\max\{|Q|,|R|\}|Q\alp-R|\ll_{\alp,\beta}1,\qquad\log\max\{|Q|,|R|\}|Q\beta-R|_p\ll_{\alp,\beta}1.\eqno\hbox{\rm(1.4)}$$ With $\beta=0$, this result means that $1/\alp$ satisfies the mixed conjecture, and was proved in \cite{dMT}.

Further questions can be considered. 
For instance, given $(\alp,\beta)\in\R\times\Q_p\backslash\{(0,0)\}$ and $\eps>0$, 
does there exist a triple of nonzero integers $(q,r,s)$ such that
 $$|qs||q\alp-r||q\beta-s|_p<\eps\ ? $$
 This condition is weaker than the multiplicative conjecture
 $$\inf_{qs\ne0}|qs||q\beta-s|_p = 0; $$
 however, we do not know additional examples.\par 
 
Inhomogeneous problems could be studied as well. 
For example, for $\eps > 0$ and a triple $(\alp_1,\alp_2,\beta)\ne(0,0,0)$ with $\alp_1, \alp_2$ in $\R$ and $\beta$ in $\Q_p$, 
does there exist a pair $(q,r)$ of nonzero integers such that 
 $$\max\{|q|,|r|\}|q\alp_1-\alp_2-r||q\beta-r|_p<\eps\ ?$$

 \section{Proof of Theorem 1.1}  
Let $f(X)$ in $\Z[X]$ be the minimal defining polynomial of $\beta$.    
We shall use units and $p$-units in the quadratic field $\Q(\beta)$.
Let us recall that a nonzero $x$ in $\Q(\beta)$ is a $p$-integer if there exists an integer $\nu$ such that $p^\nu x$ 
is an algebraic integer in $\Q(\beta)$. A nonzero $x$ in $\Q(\beta)$ is a $p$-unit if $x$ and $1/x$ are $p$-integers. Recall that 
$x$ in $\Q(\beta)$ is an algebraic integer if and only if we have $|x|_v\le1$ for each non archimedean absolute value $| \cdot |_v$ on $\Q(\beta)$, hence $x$ is a $p$-integer if and only if $|x|_v\le1$ for each non archimedean absolute value $| \cdot |_v$ on $\Q(\beta)$ such that $|p|_v=1$. Consequently, $x$ is an ordinary unit if and only if $|x|_v=1$ for each non archimedean absolute value $| \cdot |_v$ on $\Q(\beta)$, and $x$ is a $p$-unit if and only if $|x|_v=1$ for each non archimedean absolute value $| \cdot |_v$ on $\Q(\beta)$ with $|p|_v=1$.\par 

Let $U$, respectively, $U_p$, denote the multiplicative group of (ordinary) units, respectively, $p$-units, in $\Q(\beta)$. 
Let $\rho_0$, respectively, $\rho_p$, denote the number of irreducible factors of $f(X)$ in $\R[X]$, respectively, in $\Q_p[X]$. 
The group $U$ has rank $\rho_0-1$ and $U_p$ has rank $\rho_0+\rho_p-1$ (see \cite{Lang}, p. 66). Here $\rho_p=2$, since $\beta\in\Q_p$. Thus, if the field $\Q(\beta)$ has no real embedding, then $U$ has rank $0$ and $U_p$ has rank $2$; if $\Q(\beta)$ has a real embedding, then $U$ has rank $1$ and $U_p$ has rank $3$.
 \par 
 Let ${\rm id}$ and $\sig$ denote the automorphisms of $\Q(\beta)$, viewed as a subfiefd of $\Q_p$.   
 Let $\psi$ denote a complex embedding  of the field $\Q(\beta)$. Let $\tau$ be the automorphism of $\Q(\psi(\beta))$
 other than ${\rm id}$, i.e. $\tau=\psi\circ\sig\circ\psi^{-1}$. We use a method introduced by  
 Peck \cite{Pe} who obtained an effective version of a celebrated result of Cassels and Swinnerton-Dyer \cite{CS} 
 asserting that every pair in a real cubic number field  
 satisfies the Littlewood conjecture. Later, Teuli\'e \cite{Teu} used Peck's method for $p$-adic numbers. 
 Our approach is similar.  
 
  Let $\gamma$ be a nonzero number in ${\Q}(\beta)$ such that $\gamma$ and $\gamma\beta$ are algebraic integers. 
 Let $\zeta$ be a $p$-unit in $\Q(\beta)$ which also is an algebraic integer.
 We consider integers $q$ and $r$ given by the traces 
 $$q={\rm Tr}(\gamma\zeta),\qquad r={\rm Tr}(\beta\gamma\zeta).$$ 
 Note that if we consider a $p$-adic number $x\in{\Q}(\beta)$, the rational number ${\rm Tr}(x)$ can be calculated in ${\Q}_p$ as well as in $\C$, namely 
$${\rm Tr}(x)=x+\sig(x)=\psi(x)+\tau(\psi(x)).$$
 Thus, in order to estimate $q$, $r$, $|q\alp-r|$ and $|q\beta-r|_p$, we have 
 $$q=\psi(\gamma\zeta)+ \tau\circ\psi(\gamma\zeta),\qquad r=\psi(\beta\gamma\zeta)+ \tau\circ\psi(\beta\gamma\zeta),\eqno\hbox{(2.1)}$$ and also
 $$q=\gamma\zeta+ \sig(\gamma\zeta),\qquad r=\beta\gamma\zeta+ \sig(\beta\gamma\zeta).\eqno\hbox{(2.1')}$$ As $\tau(\psi(\beta))\ne\psi(\beta)$, we deduce from (2.1) that $$\max\{|q|,|r|\}\asymp_\beta\max\{|\psi(\gam\zeta) |,|\tau(\psi(\gam\zeta))|\}, \eqno\hbox{(2.2)}$$ 
 where the constants involved only depend upon $\beta$. We have
 $$q\alp-r=\psi(\gamma\zeta)(\alp-\psi(\beta)) +\tau(\psi(\gam\zeta))(\alp-\tau(\psi(\beta))),\eqno\hbox{(2.3)}$$ and by (2.1'), $$q\beta-r=(\beta-\sig(\beta))\sig(\gam)\sig(\zeta),$$ hence
$$|q\beta-r|_p\ll_\beta|\sig(\zeta)|_p, \eqno\hbox{(2.4)}$$ since $|\sig(\gam)|_p\le1$.\par
We examine successively the case where $\psi(\beta)$ is an imaginary quadratic number, and the case where $\psi(\beta)$ is a real number.

\subsection{The case of an imaginary embedding.} 

In this case, any unit in $\Q(\beta)$ has finite order. The group $U_p$ of the $p$-units of ${\Q}(\beta)$ has rank 2 (see \cite{Lang}, p. 66), hence there exists a $p$-unit $\eta$ such that $p$ and $\eta$ are multiplicatively independent. Therefore $\eta$ is not a rational number. 
Setting $\omega=\eta^2/N_{{\Q}(\beta)/{\Q}}(\eta)=\eta/\sig(\eta)$,   
we have 
$$
N_{{\Q}(\beta)/{\Q}}(\omega)=1=N_{{\Q}(\psi(\beta))/{\Q}}(\psi(\omega))=|\psi(\omega)|^2, 
$$
hence, $$|\psi(\omega)|=1.$$ But we have
$|\omega|_p\ne1$, otherwise we would also have $|\sig(\omega)|_p=1$, and then, as the non archimedean absolute values $v$ on ${\Q}(\beta)$ such that $|p|_v<1$ are $v_1=| \cdot |_p$ and $v_2=v_1\circ\sig$,  we would deduce that $\omega$ must be an ordinary unit;   
this is impossible since there are only units of finite order in ${\Q}(\beta)$. Thus, replacing if necessary $\omega$ by $\omega^{-1}$, we can suppose that $$|\omega|_p>1.$$ Let $t$ be the positive integer such that
 $$|\omega|_p=p^t,\qquad |\sig(\omega)|_p=p^{-t}.$$ Set
$p^t\omega=\omega_1$, so that
$$|\psi(\omega_1)|=p^t,\qquad|\omega_1|_p=1,\qquad |\sig(\omega_1)|_p=p^{-2t},\qquad N_{{\Q}(\beta)/{\Q}}(\omega_1)=p^{2t}.$$ Hence $\omega_1$ is a $p$-unit, and, as for each non archimedean absolute value $v$ on ${\Q}(\beta)$ such that $|p|_v<1$ we have $|\omega_1|_v\le1$, we get that $\omega_1$ is an algebraic integer.\par
For any positive integer $n$, set 
$$q_n={\rm Tr}(\gam\omega_1^n),\qquad r_n={\rm Tr}(\beta\gam\omega_1^n).\eqno\hbox{(2.5)}$$
These numbers are integers. Using (2.2), (2.3), and (2.4), we deduce from (2.5) that  
$$\max\{|q_n|,|r_n|\}\asymp_\beta|\psi(\gam)| |\psi(\omega_1^n)| =|\psi(\gam)|p^{nt} \eqno\hbox{(2.6)}$$   
and
$$|q_n\beta-r_n|_p\ll_\beta|\sig(\omega_1^n)|_p=p^{-2nt}.\eqno\hbox{(2.7)}$$ 
Put $$\psi(\omega_1)=p^te^{i\pi\theta},$$ where $\theta$ is a real number. The number $\theta$ is irrational because $p$ and $\omega_1$ are multiplicatively independent. 
For the choice of $\gam$, we simply require that $\gam$ is a positive integer (depending only upon $\beta$) such that $\gam\beta$ is an algebraic integer. Formula (2.3) becomes
$$q_n\alp-r_n=\gam p^{nt}(e^{i\pi n\theta}(\alp-\psi(\beta))+e^{-i\pi n\theta}(\alp-\overline{\psi(\beta)})),$$ 
where the bar denotes the complex conjugate.     
As $\alp\ne\psi(\beta)$, since $\alp$ lies in $\R$ while $\psi(\beta)$ is an imaginary complex number, we can write
$$|q_n\alp-r_n|\ll_{\alp,\beta} p^{nt}\left\vert e^{2i\pi n\theta}+{\alp-\overline{\psi(\beta)}\over\alp-\psi(\beta)}\right\vert.\eqno\hbox{(2.8)}$$
Let $\xi=\xi(\alp)$ denote a real number such that
$${\alp-\overline{\psi(\beta)}\over\alp-\psi(\beta)}=-e^{2i\pi\xi}.$$ 
Then, it follows from (2.8) that $$|q_n\alp-r_n|\ll_{\alp,\beta} p^{nt}\Vert n\theta-\xi\Vert,\eqno\hbox{(2.9)}$$
where $\Vert.\Vert$ is the distance to the nearest integer.  As $\theta$ is irrational, the sequence $(n\theta)_{n \ge 1}$    
is everywhere dense modulo $1$. 
Consequently, for every $\eps>0$, there are infinitely many positive integers $n$ such that $$\Vert n\theta-\xi\Vert\le\eps,\eqno\hbox{(2.10)}$$ hence by (2.9),
$$|q_n\alp-r_n|\ll_{\alp,\beta}\eps p^{nt}.\eqno\hbox{(2.11)}$$ By (2.6), (2.7), and (2.11), we have for such integers $n$ 
: $$|q_n\alp-r_n|\ll_{\alp,\beta}\eps\max\{|q_n|,|r_n|\},\qquad|q_n\beta-r_n|_p\ll_{\alp,\beta}(\max\{|q_n|,|r_n|\})^{-2},$$ hence
$$\max\{|q_n|,|r_n|\}|q_n\alp-r_n||q_n\beta-r_n|_p\ll_{\alp,\beta}\eps.$$ Thus the pair $(\alp,\beta)$ satisfies (EK+). 
A quantitative inequality can be proved by noticing that a result of Khintchine \cite{Kh} (see also \cite{JWS} or \cite{Des}) shows that for any real number $x$, and for any real number $c>5^{-1/2}$, there are infinitely many positive integers $n$ with $$\Vert n\theta-x\Vert\le{c\over n}\cdot\eqno\hbox{(2.10')}$$
By using (2.10') in place of (2.10), we get   
that there is an infinite set of positive integers $n$ which satisfy $$|q_n\alp-r_n|\ll_{\alp,\beta}{p^{nt}\over n}\cdot\eqno\hbox{(2.11')}$$
We then deduce from (2.6) and (2.11')  
that there are infinitely many positive integers $n$ for which
$$|q_n\alp-r_n|\ll_{\alp,\beta}\max\{|q_n|,|r_n|\}(\log\max\{|q_n|,|r_n|\})^{-1}, \quad |q_n\beta-r_n|_p\ll_\beta(\max\{|q_n|,|r_n|\})^{-2}.$$
Hence we have established (1.3), and Theorem 1.1 is proved in this case.

\subsection{The case of a real embedding.} 

We now suppose that the field $\Q(\beta)$ has a real embedding $\psi$. We denote by $\tau$ the automorphism of $\Q(\psi(\beta))$   
other than identity, that is, $\tau=\psi \circ\sig\circ\psi^{-1}$ where $\sig$ is the automorphism of $\Q(\beta)$ other than identity.  
\par
The group $U$ of units of ${\Q}(\beta)$ has rank 1, and the group $U_p$ of $p$-units has rank $3$ (see \cite{Lang}, p. 66). As above, we have to choose suitable units.   
\begin{Lemma}  There exist a unit $\omega_1$ and a $p$-unit $\omega_2$ with the following properties:  
$$ \psi(\omega_1)>1,\qquad 0<\psi(\sigma(\omega_1))<1,\qquad N_{{\Q}(\beta)/{\Q}}(\omega_1)=1, \eqno\hbox{\rm(2.12)}$$
$$\psi(\omega_2)>0,
\qquad|\omega_2|_p=1,\qquad|\sig(\omega_2)|_p<1,\qquad 
N_{{\Q}(\beta)/{\Q}}(\omega_2)=|\sig(\omega_2)|^{-1}_p \in p^{2 {\mathbb N}}. 
\eqno\hbox{\rm(2.13)}$$    
Moreover $\omega_2$ is an algebraic integer.
\end{Lemma}
\begin{proof}
First we take a unit $\widehat\omega_1$ in $U$ of infinite order and set $\omega_1 = \widehat\omega_1^2$.
Then, we have $$|\omega_1|_p=|\sig(\omega_1)|_p=1,\qquad\psi(\omega_1)>0,\qquad N_{\Q(\beta)/\Q}(\omega_1)=1.$$ Replacing if necessary $\omega_1$ by $\omega_1^{-1}$, we get (2.12). Now, as $U_p$ has rank $3$, there exists a $p$-unit $\widehat\omega_2$ such that $p$, $\omega_1$ and $\widehat\omega_2$ are multiplicatively independent. Denote by $u$ the integer such that $|\widehat\omega_2|_p=p^{u}$, and put
  $p^{2u}\widehat\omega_2^2 =\omega_2$, so that $|\omega_2|_p=1$. Then
 $\omega_2 $ is a $p$-unit. Since $\omega_1$ and $\omega_2$ are multiplicatively independent, and since $U$ has rank one, 
 $\omega_2$ cannot be an ordinary unit. Therefore we must have $|\sig(\omega_2)|_p\ne1$.
 otherwise we would have $|\omega_2|_v=1$ for any non-archimedean absolute value $| \cdot |_v$, and thus $\omega_2$ would be an ordinary unit. 
Replacing if necessary $\omega_2$ by $\omega^{-1}_2$, we can suppose that $|\sig(\omega_2)|_p<1$.
Since $N_{{\Q}(\beta)/{\Q}}(\omega_2) = (N_{{\Q}(\beta)/{\Q}}(p^u \widehat\omega_2))^2$ belongs to $p^{2\Z}$, we have $$N_{{\Q}(\beta)/{\Q}}(\omega_2)=|N_{{\Q}(\beta)/{\Q}}(\omega_2)|^{-1}_p=|\sig(\omega_2)|^{-1}_p,\eqno\hbox{(2.14)}$$ and thus we get (2.13). Moreover, $\omega_2$ is an algebraic integer, since $|\omega_2|_v\le1$ for each non archimedean absolute value on ${\Q}(\beta)$.
\end{proof} 

As above, let $\gamma$ be a nonzero element in ${\Q}(\beta)$ such that $\gamma$ and $\gamma\beta$ are algebraic integers. We shall see that we must choose $\gam=\gam(\alp,\beta)$ such that moreover $\sig(\gam)=(-1)^i\gam$ where $i=i(\alp,\beta)$ is an integer such that $$(-1)^i(\alp-\psi(\beta))(\alp-\tau(\psi(\beta)))<0.\eqno\hbox{(2.15)}$$ This is possible since, by assumption, $\alp$ is not a root of the minimal defining polynomial of $\beta$,  
therefore $\alp\ne\psi(\beta)$ and $\alp\ne\tau(\psi(\beta))$. 
Thus, $\gam$ is a nonzero integer or an algebraic integer such that $\gam^2$ lies in $\Z$ and $\gam$ is not in $\Z$. 
\par
Let $\omega_1$ be a unit and $\omega_2$ a $p$-unit such that the conditions of Lemma 2.1 are satisfied. Write 
$$|\sig(\omega_2)|_p=p^{-2t},$$ where $t$ is a positive integer. For a pair $(n_1,n_2)$ of integers with $n_2 > 0$, consider 
 $$q_{n_1,n_2}={\rm Tr}(\gam\omega_1^{n_1}\omega_2^{n_2}),\qquad  r_{n_1,n_2}={\rm Tr}(\beta\gam\omega_1^{n_1}\omega_2^{n_2}).$$ These formulae are identical to (2.1) with $\zeta=\omega_1^{n_1}\omega_2^{n_2}$.
Since $\omega_1$ is a unit and $\omega_2$ an algebraic integer, the rational numbers $q_{n_1,n_2}$ and $r_{n_1,n_2}$ are integers,  
and by (2.2) and (2.4), we have
$$\max\{|q_{n_1,n_2}|,|r_{n_1,n_2}|\}\asymp_\beta\max\{|\psi(\omega_1^{n_1}\omega_2^{n_2})|,|\tau(\psi(\omega_1^{n_1}\omega_2^{n_2}))|\}\eqno\hbox{(2.16)}$$ and
$$|q_{n_1,n_2}\beta-r_{n_1,n_2}|_p\ll_\beta|\sig(\omega_1^{n_1}\omega_2^{n_2})|_p=p^{-2n_2t}.\eqno\hbox{(2.17)}$$
By (2.3), we also have
$$q_{n_1,n_2}\alp-r_{n_1,n_2}=\psi(\gam) \bigl(\psi(\omega_1^{n_1}\omega_2^{n_2})(\alp-\psi(\beta))+(-1)^i\tau(\psi(\omega_1^{n_1}\omega_2^{n_2}))(\alp-\tau(\psi(\beta))) \bigr),$$
and by (2.15), we get
$$|q_{n_1,n_2}\alp-r_{n_1,n_2}|\ll_\beta\bigl\vert\psi(\omega_1^{n_1}\omega_2^{n_2})|\alp-\psi(\beta)|-\tau(\psi(\omega_1^{n_1}\omega_2^{n_2}))|\alp-\tau(\psi(\beta))|\bigr\vert. \eqno\hbox{(2.18)}$$    
 Recalling that by (2.12), (2.13), and (2.14), we have $N_{\Q(\beta)/\Q}(\omega_1^{n_1}\omega_2^{n_2})=p^{2n_2t}$, we rewrite (2.18) in the form
 $$|q_{n_1,n_2}\alp-r_{n_1,n_2}|\ll_{\alp,\beta}{p^{2n_2t}\over\psi(\omega_1^{n_1}\omega_2^{n_2})}\left|
\psi(\omega_1^{2n_1}\omega_2^{2n_2})p^{-2n_2t}-{|\alp-\tau(\psi(\beta)|\over|\alp-\psi(\beta)|}\right|.\eqno\hbox{(2.18')}$$
The real numbers $\log\psi(\omega_1)$ and $\log\psi(\omega_2)-t\log p$ are linearly independent over $\Z$ since $\omega_1$, $\omega_2$ and $p$ are multiplicatively independent. Thus the sequence    
 $(n(\log\psi(\omega_2)-t\log p)/\log\psi(\omega_1))_{n \ge 1}$ 
 is everywhere dense modulo $1$. Consequently, for every real number $\eps$ with $0<\eps\le1$, there exists a pair $(n_1,n_2)$ of integers, with $n_2>0$, such that
$$|2n_1\log\psi(\omega_1)+2n_2(\log\psi(\omega_2)-t\log p)-\log\vert\alp-\tau(\psi(\beta))|+\log|\alp-\psi(\beta)|\vert<\eps.\eqno\hbox{(2.19)}$$ In particular, we have
$$-{1\over2}\le\log{\psi(\omega_1^{n_1}\omega_2^{n_2})\over p^{n_2t}} -\log\left\vert{\alp-\tau(\psi(\beta))\over\alp-\psi(\beta)}\right\vert^{1/2}\le{1\over2},$$  hence
$$e^{-1/2}\left\vert{\alp-\tau(\psi(\beta))\over\alp-\psi(\beta)}\right\vert^{1/2}\le \psi(\omega_1^{n_1}\omega_2^{n_2})p^{-n_2t}\le e^{1/2}\left\vert{\alp-\tau(\psi(\beta))\over\alp-\psi(\beta)}\right\vert^{1/2}.\eqno\hbox{(2.20)}$$
Thus, by (2.18') and (2.20), we have
$$|q_{n_1,n_2}\alp-r_{n_1,n_2}|\ll_{\alp,\beta}p^{n_2t}
\left|\psi(\omega_1^{2n_1}\omega_2^{2n_2})p^{-2n_2t}-{|\alp-\tau(\psi(\beta))|\over|\alp-\psi(\beta)|}\right|,\eqno\hbox{(2.21)}$$
and, by using again (2.20), we deduce from (2.19) and (2.21) that
$$|q_{n_1,n_2}\alp-r_{n_1,n_2}|\ll_{\alp,\beta}\eps p^{n_2t}.\eqno\hbox{(2.22)}$$
Since $N_{\Q(\beta)/\Q}(\omega_1^{n_1}\omega_2^{n_2})=p^{2n_2t}$, we deduce from (2.20) that we also have $$e^{-1/2}\left\vert{\alp-\psi(\beta)\over\alp-\tau(\psi(\beta))}\right\vert^{1/2}\le \tau(\psi(\omega_1^{n_1}\omega_2^{n_2}))p^{-n_2t}\le e^{1/2}\left\vert{\alp-\psi(\beta)\over\alp-\tau(\psi(\beta))}\right\vert^{1/2},\eqno\hbox{(2.20')}$$ hence, by (2.16), (2.17), (2.20), and (2.20'), we get 
$$\max\{|q_{n_1,n_2}|,|r_{n_1,n_2}|\}\asymp_{\alp,\beta}
 p^{n_2t},\qquad |q_{n_1,n_2}\beta-r_{n_1,n_2}|_p\ll_\beta p^{-2n_2t} .\eqno\hbox{(2.23)}$$
Thus,
(2.22) and (2.23) ensure that
 $$ \max\{|q_{n_1,n_2}|,r_{n_1,n_2}|\}|q_{n_1,n_2}\alp-r_{n_1,n_2}||q_{n_1,n_2}\beta-r_{n_1,n_2}|_p\ll_{\alp,\beta}\eps,$$ which proves that 
 the pair $(\alp, \beta)$ satisfies (EK+).    
\par As above, the Theorem of Khintchine \cite{Kh} shows that inequality (2.19) holds true for 
$\eps\ll_\beta 1/n_2$ for pairs of integers $(n_1,n_2)$ with $n_2$ arbitrarily large. Since, by (2.23),
$$n_2\asymp_{\alp,\beta}\log\max\{|q_{n_1,n_2}|,|r_{n_1,n_2}|\},$$ inequalities (2.22) and (2.23) lead to (1.3), which proves Theorem 1.1.

\section{Proof of Theorem 1.2.}
Denote by $\tau$ the automorphism of the field $\Q(\alp)$ other than identity.   
Let $\vph$ be an isomorphism of $\Q(\alp)$ into the algebraic closure $\Omega_p$ of $\Q_p$. Also, denote by $\sig$ the automorphism 
of $\Q(\vph(\alp))$ other than identity, that is to say $\sig=\vph\circ\tau\circ\vph^{-1}$. Notice that in the case where $\vph(\alp)$ does not lie in $\Q_p$, there exists an isometrical automorphism $\sig_p$ of $\Q_p(\vph(\alp))$ such that $\sig_p(\vph(\alp))=\sig(\vph(\alp))$, that is to say that $\sig_p(x)=\sig(x)$ for each $x\in\Q(\vph(\alp))$. We will set again $\sig_p=\sig$.\par
Let $\zeta$ be a unit of infinite order in $\Q(\alp)$. 
Replacing if necessary $\zeta$ by $\zeta^{\pm2}$ we can suppose that $\zeta>1$ and $0<\tau(\zeta)=\zeta^{-1}<1$.
As $|\vph(\zeta)|_p=1$, replacing again $\zeta$ by $\zeta^\nu$, where $\nu$ is a suitable  
positive integer, we may suppose that $$|\vph(\zeta)-1|_p<p^{-1/(p-1)}.\eqno\hbox{(3.1)}$$ Indeed, noticing that $|\vph(\zeta^{p^2-1})
-1|_p\le p^{-1/2}$, it is easy to see that we can take $\nu=p^2-1$ if $p\ge5$, $\nu=24$ if $p=3$, and $\nu=12$ if $p=2$. 
\par Let $\gamma$ be in $\Q(\alp)$ such that $\gam$ and $\gam\alp$ are algebraic integers, and $\gam>0$. For any positive integer $n$, we set
$$q_n={\rm Tr}(\gam\zeta^n)=\gam\zeta^n+\tau(\gam)\zeta^{-n},\qquad r_n={\rm Tr}(\alp\gam\zeta^n)=\alp\gam\zeta^n+\tau(\alp\gam)\zeta^{-n}.$$ These numbers are integers, and, for $n$ sufficiently large (depending on $\gam$), we have $q_n>0$ and $$q_n\asymp_{\alp,\gam}\zeta^n. \eqno\hbox{(3.2)}$$ As
$$q_n\alp-r_n=(\alp-\tau(\alp))\tau(\gam)\zeta^{-n},$$ we also have $$|q_n\alp-r_n|\asymp_{\alp,\gam}\zeta^{-n}.$$ Thus
$$q_n|q_n\alp-r_n|\asymp_{\alp,\gam}1.\eqno\hbox{(3.3)}$$
Let us do similar computations in $\Omega_p$ and write     
 $$q_n={\rm Tr}(\varphi(\gam\zeta^n))=\vph(\gam\zeta^n)+\sig(\vph(\gam\zeta^n)), \qquad r_n={\rm Tr}(\varphi(\alp\gam\zeta^n))=\vph(\alp\gam\zeta^n)+\sig(\vph(\alp\gam\zeta^n)).$$ 
 As $\sig(\vph(\zeta))=\vph(\tau(\zeta))=\vph(1/\zeta)$, we get 
$$q_n\beta-r_n=(\beta-\vph(\alp))\vph(\gam)\vph(\zeta)^n+(\beta-\sig(\vph(\alp)))\sig(\vph(\gam))\vph(\zeta)^{-n}.$$ 
By assumption, $\beta$ is not a root of the minimal defining polynomial of $\alpha$, thus  
 $\beta\ne\vph(\alp)$ and $\beta\ne\sig(\vph(\alp))$. Since $|\vph(\zeta)|_p=1$, we can write
$$|q_n\beta-r_n|_p\asymp_{\alp,\beta,\gam}\left|\vph(\zeta)^{2n}+{(\beta-\sig(\vph(\alp)))\sig(\vph(\gam))\over(\beta-\vph(\alp))\vph(\gam)} \right|_p.\eqno\hbox{(3.4)}$$ 
 We have to choose a suitable $\gam$.  
 
 \begin{Lemma} There exists a nonzero number $\gam\in\Q(\alp)$ such that
 $$\left\vert{(\beta-\sig(\vph(\alp)))\sig(\vph(\gam))\over(\beta-\vph(\alp))\vph(\gam)}+1\right\vert_p<|\vph(\zeta)-1|_p.\eqno\hbox{\rm(3.5)}$$
 Moreover, we can suppose that $\gam$ and $\alp\gam$ are algebraic integers, and $\gam>0$.
 \end{Lemma}
 
\begin{proof}
First, note that it is enough to find a nonzero $\gam$ in $\Q(\alp)$ which satisfies (3.5), since $\sig(\gam)/\gam$ is unchanged if we replace $\gam$ by $G\gam$, where $G$ is a nonzero integer.
This allows us to suppose that $\gam$ is an algebraic integer. Notice also that, if $\beta-\vph(\alp)=-(\beta-\sig(\vph(\alp)))$, i.e., if $\beta={\rm Tr}(\alp/2)$, we have the obvious solution $\gam=1$. Thus we can suppose that $\beta\ne{\rm Tr}(\alp/2)$.\par
Let us consider the homographic function $H$ on $\Q_p(\vph(\alp))$ given by    
 $$H(t)={t+\sig(\vph(\alp))\over t+\vph(\alp)}\cdot$$ Since $\sig(\vph(\alp))\ne\vph(\alp)$, this function is not constant and takes 
 in $\Q_p(\vph(\alp))$ every value different from $1$.   
Hence there exists $t_1\in\Q_p(\alp)$ such that $$H(t_1)=-{\beta-\vph(\alp)\over\beta-\sig(\vph(\alp))}\cdot$$ 
We check that 
$$t_1={-\beta{\rm Tr}(\alp)+{\rm Tr}(\alp^2)\over2\beta-{\rm Tr}(\alp)}\cdot$$
Notice that $t_1$ lies in $\Q_p$. Consequently, as the function $H$ is continuous at the point $t_1$, and since $\Q$ is everywhere dense in $\Q_p$, there exist integers $x$ and $y$ with $y\ne0$, such that
$$\left|H\left({x\over y}\right)+{\beta-\vph(\alp)\over\beta-\sig(\vph(\alp))}\right|_p<\left|
{\beta-\vph(\alp)\over\beta-\sig(\vph(\alp))}\right|_p|\vph(\zeta)-1|_p.$$
This means that, if we put $\gam=x+y\alp$, we have $\gam\ne0$ and
$$\left|{\sig(\vph(\gam))\over\vph(\gam)}+{\beta-\vph(\alp)\over\beta-\sig(\vph(\alp))}\right|_p<\left|{\beta-\vph(\alp)\over\beta-\sig(\vph(\alp))}\right|_p|\vph(\zeta)-1|_p,$$
which is (3.5).
\end{proof}

We can now complete the proof of Theorem 1.2. By Lemma 3.1, we choose $\gam > 0$ satisfying  
(3.5) and such that $\gam$ and $\alp\gam$ are algebraic integers. 
We use the $p$-adic logarithm function $\log_p$. 
For
 $z\in\Q_p(\vph(\alp))$ with $|z|_p<p^{-1/(p-1)}$, we have
 $$\log_p(1+z)=\sum_{n=1}^{+\infty}{(-1)^{n-1}z^n\over n},$$ (see \cite{Rb}).
 The function $z\longmapsto\log_p(1+z)$ is a map from the ball $\{ z \in  \Q_p(\vph(\alp)) : |z|_p<p^{-1/(p-1)} \}$
 into itself which is isometrical. In view of the choices of $\zeta$ and $\gam$, putting 
 $$\lambda=-{(\beta-\sig(\vph(\alp))) \sig(\vph(\gam)) \over(\beta-\vph(\alp))\vph(\gam)},$$ we have $$|\lambda-1|_p<p^{-1/(p-1)}. 
 $$
This and (3.1) show that $\log_p\lambda$ and $\log_p\vph(\zeta)$ are well defined.  
We have
 $$|\vph(\zeta)^{2n}-\lambda|_p=|2n\log_p\vph(\zeta)-\log_p\lambda|_p=|2\log_p\vph(\zeta)|_p\left|n-{\log_p\lambda\over2\log_p\zeta}\right|_p\cdot\eqno\hbox{(3.6)}$$
 Let us prove that $(\log_p\lambda)/\log_p\vph(\zeta)$ 
 belongs to $\Q_p$. This is obvious  
 if $\vph(\alp)\in\Q_p$. If $\vph(\alp)\not\in\Q_p$, then $\sig$ can be regarded as a $\Q_p$-automorphism of $\Q_p(\vph(\alp))$ and this automorphism is isometrical. We have $\sig(\vph(\zeta))=1/\vph(\zeta)$, and as $\sig\circ\sig={\rm id}$, we also have $\sig(\lambda)=1/\lambda$, since $\beta$ lies in $\Q_p$. Now for each $Z\in\Q_p(\vph(\alp))$ with $|Z-1|_p<p^{-1/(p-1)}$, we have $\sig(\log_pZ)=\log_p(\sig(Z))$ and $\log_p(1/Z)=-\log_pZ$. Consequently, we have 
 $\sig(\log_p\lambda)=-\log_p\lambda$ and  $\sig(\log_p\vph(\zeta))=-\log_p(\vph(\zeta))$. Thus we get
 $$\sig\left({\log_p\lambda\over\log_p\vph(\zeta)}\right)={\log_p\lambda\over\log_p\vph(\zeta)},$$ which proves that 
${\log_p\lambda\over\log_p\vph(\zeta)}$ lies in $\Q_p$. Moreover, as $|{\log_p\lambda\over\log_p\vph(\zeta)}|_p={|\lambda-1|_p\over|\vph(\zeta)-1|_p}<1$, by (3.5), the quotient  
${\log_p\lambda\over2\log_p\vph(\zeta)}$ belongs to $\Z_p$. Therefore, for each integer $N>0$, there exists an integer $n$, with $0\le n<p^N$, such that
 $$\left|n-{\log_p\lambda\over2\log_p\vph(\zeta)}\right|_p\le p^{-N}.\eqno\hbox{(3.7)}$$    
 Thus, it follows from (3.3), (3.4), (3.6), (3.7), and (3.2) that 
 $$\max\{|q_n|,|r_n|\}|q_n\alp-r_n|\asymp_{\alp,\beta}1,\qquad |q_n \beta 
 -r_n|_p\le p^{-N},\qquad \log\max\{|q_n|,|r_n|\}\ll_{\alp,\beta} p^N,$$ 
which leads to (1.4) 
 and concludes the proof.

\bigskip

\noindent bugeaud@math.unistra.fr \hfill bernard.demathan@gmail.com

\medskip 

\noindent I.R.M.A., UMR 7501

\noindent Universit\'e de Strasbourg et CNRS  

\noindent 7 rue Ren\'e Descartes  

\noindent  67084 Strasbourg Cedex, France

\medskip 

\noindent  Institut universitaire de France

 \end{document}